\input amstex
\documentstyle{amsppt}
\input bull-ppt


\NoBlackBoxes

\topmatter
\cvol{31}
\cvolyear{1994}
\cmonth{October}
\cyear{1994}
\cvolno{2}
\cpgs{216-222}
\title A GENERAL DECOMPOSITION THEORY FOR RANDOM CASCADES  
\endtitle
\author Edward C. Waymire and  Stanley C. Williams 
\endauthor
\shortauthor{E. C. Waymire and S. C. Williams}
\address Mathematics Department, Oregon State University, 
Corvallis,  Oregon
97330\endaddress
\ml waymire\@math.orst.edu\endml
\address Mathematics Department, Utah  State University, 
Logan, Utah 84321\endaddress
\ml williams\@sunfs.math.usu.edu\endml
\date June 14, 1993\enddate
\abstract This announcement describes a probabilistic 
approach to
cascades which, in addition to providing an entirely  
probabilistic proof 
 of the Kahane-Peyri\`ere theorem for independent cascades,
 readily applies to
general dependent cascades. 
 Moreover, this unifies various seemingly disparate cascade
decompositions, including  Kahane's T-martingale 
decomposition 
and  dimension disintegration.\endabstract
\endtopmatter
\document

\heading 1.  Brief history of the problem \endheading

 A theory of 
{\it positive T-martingales} was introduced in \cite{K3} 
as the
general framework for independent multiplicative cascades 
and
random coverings.  This includes spatial distributions of 
interest
in both probability theory and in the physical sciences, 
e.g. \cite{CCD,
CK, {DE}, {DG},
{DF}, {DM}, {F}, {GW}, {MS}, {MW},
{TLS}, {PW}, S}.  For
basic definitions, let $T$ be a locally compact metric 
space with Borel sigmafield
${\Cal B,}$ and let $(\Omega, {\Cal F}, P)$ be a  
probability space together
with an increasing sequence ${\Cal F}_n, n = 1,2,\dots$, 
of sub-sigmafields of
${\Cal F}.$   
 A {\it positive T-martingale} is
a sequence $\{Q_n\}$ of non-negative functions on $T\times 
\Omega$ such
that: (i) for each $t\in T, \{Q_n(t,\cdot): n = 
1,2,\dots\}$ is a martingale
adapted to ${\Cal F}_n, n = 1,2,\dots;$ (ii) for $P$-a.s.\ 
$\omega \in \Omega,
\{Q_n(\cdot,\omega): n = 1,2,\dots\}$ is a sequence of 
Borel measurable 
nonnegative real-valued functions on $T$.

Let $M^+(T)$ denote the space of positive Borel measures
on $T$, and suppose that $\{Q_n(t)\}$ is a positive 
$T$-martingale.  For
$\sigma\in M^+(T)$ let $Q_n\sigma$ denote the random 
measure defined
by 
$Q_n\sigma << \sigma$ and ${dQ_n\sigma\over d\sigma}(t) := 
 Q_n(t), t\in T.$
Denote the space of bounded continuous functions on $T$ by 
$C_B(T).$
Then, essentially by the martingale convergence theorem,
one obtains  a random Borel measure
 $Q_\infty\sigma$ such that with probability one \cite{K3}
$$\lim_{n\to \infty}\int_T f(t)Q_n\sigma\,(dt) = \int_T 
f(t)Q_\infty\sigma
\,(dt),\qquad f\in C_B(T). \tag1.1$$
As suggested by the notation, one may view
 $\sigma\to \sigma_\infty \equiv Q_\infty\sigma$ as a 
random operator
acting on $M^+(T).$  
The following special case of multiplicative cascades is 
central to the general
theory developed in  \cite{K3}.

\subheading{Independent cascades}
Let $b\ge 2$ be a natural number,
referred to as a {\it branching number};
 and let $T = \{0,1,\dots,b-1\}^{\bold N}$ be the metric
space for the ultrametric $\rho(s,t) = b^{-a(s,t)}$ where
$a(s,t) = \inf\{ n : s_n \neq t_n\}, s =
(s_1,s_2,\dots), t = (t_1,t_2,\dots).$  The countable set
  $T^*\equiv T^*(\infty) :=
\cup_{n=0}^\infty T^*(n),$ where $T^*(n) := 
\cup_{k=1}^n\{0,1,\dots,b-1\}^k,$
provides a convenient labelling of the vertices of the  
infinite $b$-ary tree. It is sometimes convenient to 
adjoin a root vertex
denoted $\emptyset$ to $T^*.$
 For $t\in T, n\ge 1,$
write $t|n := (t_1,t_2,\dots,t_n)\in T^*, t|0 := 
\emptyset.$ For
 $\gamma = (\gamma_1,\dots,\gamma_n) \in
T^*, j\in \{0,1,\dots,b-1\},$ define $ |\gamma| := n$ and 
$\gamma*j :=
(t_1,\dots,t_n,j).$  Also denote the  $n$\<th generation 
partition by 
 $\Delta_\gamma
:= \{t\in T: t|_{|\gamma|} = \gamma\}.$  In this context 
the {\it cascade generators}
are furnished  by a  (denumerable) family  
of i.i.d.\
nonnegative mean-one random variables $\{W_{\gamma}: 
\gamma\in T^*\}$
indexed by the tree and defined
 on a probability space
 $(\Omega, {\Cal F}, P).$ With this, let
$$P_n(t) :=
\sum_{|\gamma| = n}
 W_{\gamma}{\bold 1}_{\Delta_\gamma}(t). \tag1.2$$
Now  the {\it homogeneous independent cascade}
 is the positive $T$-martingale defined
by  $Q_n = P_1P_2\cdots P_n$
with respect to 
${\Cal F}_n := \sigma\{W_{\gamma} : |\gamma| \le n\}.$   
Define 
$\lambda_\infty : = Q_\infty\lambda,$ where here and
throughout  $\lambda$ denotes  the Haar measure on $T$.
It may happen that $\lambda_\infty = 0$ a.s., referred to
as a {\it degeneracy}.  The onset of nondegeneracy may be 
viewed
as a critical phenomenon associated with the cascade.
The fundamental theorem of
 [KP] on the structure of homogeneous  independent 
cascades provides
necessary and sufficient conditions on the distribution of 
the generators
for {\it \RM{(i)} nondegeneracy and, \RM{(ii)}
 divergence of moments of $\lambda_\infty[0,1],$}
and provides  in all cases of nondegeneracy  (iii)
{\it the  Hausdorff
dimension of the  support.} Replacement of
the generators at each generation
by i.i.d.\ nonnegative random vectors ${\bold W} := 
(W_0,\dots,W_{b-1}),$
$E{1\over b}\sum W_j = 1,$ is possible  without 
substantial changes in the
conditions \cite {N}.  Moreover, this makes it possible to 
include the
random distribution functions studied by \cite{DF}.
The results described in this note provide an entirely 
probabilistic approach to the 
computation of the fine-scale properties of random 
cascades, which also
includes a simple  new solution  to
problems  (i)--(iii) for independent cascades. 
As will be seen, the power of this 
approach for more general cascades 
is that fine-scale computations are reduced 
 to  the law of large numbers (ergodic theory) problems.

Much of the recent focus of the  study of independent 
cascades, both
theoretically and empirically, has concerned various 
characteristics of their
singularity structure, e.g.\ \cite{HW, CK, F, PW,
MS, O}.
In this regard, for an arbitrary Borel measure $\sigma,$
one has a unique disintegration $\sigma(\cdot) = \int
\sigma_\alpha(\cdot)\nu\,(d\alpha)$ (properly 
interpreted), where
$\int_{[0,\beta]}\sigma_\alpha\nu\,(d\alpha)$ is supported 
by a set
of Hausdorff dimension no larger than $\beta$,
and 
$\int_{(\beta,\infty]}\sigma_\alpha(B)\nu\,(d\alpha) > 0 
\Rightarrow
\dim (B) > \beta$ \cite {C, KK}.  The measure  $\nu$  is 
called  the
{\it dimension spectrum.}  Each spectral mode 
$\sigma_\alpha$ is
supported on a set of dimension at most $\alpha.$ If $\nu$ 
has an
atom at $\alpha$, then $\sigma_\alpha$ is {\it 
unidimensional}; i.e.\
the dimension spectrum of $\sigma_\alpha$ is a Dirac point 
mass.
It follows from  \cite{KP} that the
(homogeneous) independent cascades are a.s.\ 
unidimensional.  The richness
of 
extensions of the cascade theory to dependent cascades is 
partly 
reflected in an often  natural nontrivial spectral 
disintegration.  

\heading 2.  Main results\endheading

  The original proof of the Kahane-Peryi\`ere
theorem is  based on a combination of 
probabilistic and analytic computations which make strong 
use of
the {\it statistical independence}.  The  point of focus 
in \cite{KP} is
a {\it distributional fixed-point equation} for the total 
mass of\ the
cascade.
 Extensions  to a limited class
of dependent cascades, namely, finite-state Markov 
generators, have been found
along similar lines in  \cite{WW1}.  However, even in the 
countably infinite-state Markov
context, the fixed-point analysis does not readily extend 
to dependent 
generators.  It
is interesting to note that the fixed-point equation also 
arises in certain
interacting particle context, \cite{HL, DL}. The approach 
described in
this note may be of independent interest when applied to 
the fixed-point
 problem under dependence.

To define the general {\it dependent random cascade,}  we 
begin with
a probability measure $p_0(dx)$ and  a collection
of mean-one transition probability kernels
 $q_n(dx|x_0,\dots,x_{n-1}),$ on ${\Cal B} [0,\infty),$ 
 $x_i \ge 0, n \ge 1.$  Using the Kolmogorov extension 
theorem, one may
construct a unique probability measure $P$ on the product 
space
$\Omega \coloneq [0,\infty)^{T^*(\infty)},
{\Cal \coloneq := {\Cal B}^{T^*(\infty)}[0,\infty)}$ 
together with the coordinate
projection random variables $W_\gamma(\omega) \coloneq 
\omega_\gamma$ such
that (i) $W_\emptyset$ has distribution $p_0(dx);$ (ii) 
for any $t\in T,$ 
the conditional distribution of
$W_{t|n+1}$ given $W_\gamma, |\gamma|\le n,$ is
 $q_{n+1}(dx|W_{t|0},W_{t|1},\dots,W_{t|n});$ and (iii) 
for $s\in T,
s|n+1\neq t|n+1, W_{s|n+1}$ is conditionally independent 
of $W_{t|n+1}$
given $W_\gamma, |\gamma|\le n$ \cite{WW2}.
We refer to this model as the
{\it cascade generator} corresponding to the given 
transition kernels 
$q_n$ and
initial distribution $p_0.$
By relabeling the states,
 the two-state \lq\lq Markov chains on trees\rq\rq \ 
constructed in
 \cite{Pr, Sp} may be adapted as a special case to 
illustrate this
general setting.

Given a cascade generator, one can define 
 a positive $T$-martingale  by applying
the formula (1.2). 
The resulting random measure $\lambda_\infty := 
Q_\infty\lambda$ will be referred to
as the {\it dependent cascade.} 
Our approach to the study of cascades in this generality 
is based on three
elements, namely, (i) {\it a weight system perturbation,}  
(ii) {\it a
size-bias transform\/}, and (iii) {\it a 
general percolation method}.

\dfn {Definition 2.1} A {\it weight system} is a family 
${\bold F}$
of real-valued functions $F_{\gamma}:\Omega\to [0,\infty)$
indexed by the tree,  where
for each $\gamma\in T^*,$ $F_{\gamma}$ is  
$\sigma\{W_{\gamma|j}:
j\le |\gamma|\}$ measurable, such that
$Q_{{\bold F},n}(t) := \sum_{|\gamma|=n}\prod_{j\le 
n}(W_{\gamma|j})F_{\gamma}{\bold
1}_{\Delta_\gamma}$ is a positive $T$-martingale, referred 
to
as a {\it weighted cascade.} A weight
system ${\bold F}$ is called a {\it weight decomposition} 
in the case ${\bold
F}^c := \{1-F_\gamma:\gamma\in T^*\}$ is also a weight 
system. 
\enddfn

Note that a weight system is a weight decomposition if and 
only if
the weights are bounded between 0 and 1.  

Assuming that one already has a nondegenerate limit cascade,
Peyri\`ere \cite{P}  defines a probability ${\Cal Q}$
on $\Omega\times T$, named the {\it Peyri\`ere
probability} in \cite{K1}, for the joint distribution
of a randomly selected path and the cascade generators 
along the path. 
This probability plays an important role in the analysis 
of the structure
of independent cascades {\it provided} nondegeneracy has 
been established. 
The following probability is a useful extension of this 
notion in 
two directions: namely, (i) it
does not require an a priori nondegeneracy condition; and 
(ii) it permits
perturbations by a weight system.
 For a given weight system 
${\bold F,}$
let $Q_{{\bold F},n}(t) := b^{-n}\prod_{j\le n}( W_{t|j}) 
F_{t|n},$ $t\in T,
 n\in {\bold N}.$  Define
a sequence ${\Cal Q}_{{\bold F},n}$ of measures on 
$\Omega\times T$ by
$$\int_{\Omega\times T}f(\omega,t){\Cal Q}_{{\bold 
F},n}\,(d\omega\times dt)
:= E_P\int_Tf(\omega,t)Q_{{\bold F},n}(t)\lambda 
\,(dt),\tag2.1$$
for bounded measurable functions $f$.  Then one normalizes 
the masses
of the ${\Cal Q}_{{\bold F},n}$ by a factor $Z_\emptyset 
:= EW_{\emptyset}F_{\emptyset}$ and
extends to a probability ${\Cal Q}_{{\bold F}}$ using the 
Kolmogorov
extension theorem.  One requires $Z_\emptyset  > 0$ here 
and throughout.

The third ingredient to this theory is a generalization of 
an idea 
considered without proof in 
\cite {K1}, which may be viewed as a percolation method; 
see \cite{WW2} 
for proofs.  By
independently  pruning the tree, one studies the critical 
parameters governing 
the
survival of mass,
i.e. nondegeneracy of the percolated cascade, to determine 
dimension
estimates on the support of the cascade.  This is similar 
to an idea explored in
\cite{L} but differs by the distinction between {\it 
locations}  and {\it
amounts} of positive mass.

The following results provide the foundations for the 
general theory being
announced here.

\proclaim{ Theorem 2.1 \rm(A Lebesgue decomposition)}  Let 
$\pi_\Omega$
denote the coordinate projection map of  $\Omega\times T$ 
onto
$\Omega.$  Then
$$d{\Cal Q_{\bold 
F}}\circ\pi_\Omega^{-1}=Z_\emptyset^{-1}\lambda_{{\bold 
F},\infty}(T){\bold 1}(\lambda_{{\bold
F},\infty}(T) <  \infty)dP + 
{\bold 1}(\lambda_{{\bold F},\infty}(T) = \infty)
d{\Cal Q}_{\bold F}\circ\pi_\Omega^{-1}$$
where $\lambda_{{\bold F},\infty} = Q_{{\bold 
F},\infty}\lambda.$

\endproclaim

\proclaim{ Theorem 2.2 \rm (A size-biased disintegration)}
Given a weighting system ${\bold F},$
$${\Cal Q_{\bold F}}(d\omega,dt) = P_{{\bold 
F},t}(d\omega)\lambda(dt),$$
where 
$${dP_{{\bold F},t}\over dP}|_{{\Cal F}_n}  = \prod_{j\le 
n}(W_{t|j})F_{t|n}.$$
\endproclaim

\proclaim{Theorem 2.3 \rm(A submartingale bound)} 
Let ${\Cal F}_{t,n}
:= \sigma\{W_{t|i}, W_\gamma \!:\! i\! \ge \!0, 
|\gamma|\!\le\! n\}.$ Fix $c_k \ge 0,$
such that $c_k$ is ${\Cal F}_{t,0}$-measurable,
 $E_{{\Cal Q}_{\bold F}}\sum_{k}c_k < \infty.$
Given a weight
system ${\bold F},$ 
letting $\lambda_{{\bold F},n} = Q_{{\bold F},n}\lambda,$ 
one has 
 for arbitrary $t\in T$
$$b^{-n}\prod_{i\le n}( W_{t|i})F_{t|n} \le 
\lambda_{{\bold F},n}(T)
\le b^{-n}\prod_{i\le n}( W_{t|i})F_{t|n} +
\sup_{j}({\prod_{i\le j}(W_{t|i}b^{-1})\over c_j})M_n,$$
where 
$$M_n = \sum_{j=0}^{n-1}c_jb^{-(n-j)}\sum_{|\gamma| = 
n-j\atop \gamma(1)\neq
t(j+1)} \prod_{i=1}^{n-j}(W_{[t|j]*[\gamma|i]})F_\gamma$$
is a nonnegative $P_{{\bold F},t}$-submartingale with 
respect to 
${\Cal F}_{t,n},$
whose Doob decomposition has the predictable part
 $A_n = {b-1\over b}\sum_{j\le n-1} c_j.$
\endproclaim

\proclaim{Theorem 2.4 \rm(A percolation method)}
 Let $Q_{\beta,n}(t), t\in T,$ be
the $\beta$-model defined by the  cascade with independent
 generators 
$$ W_\gamma := B_\gamma  = 
\cases b^\beta &\text{w. prob.
 $b^{-\beta}$},\\
0 &\text{else}\endcases$$
and independent of the weighted cascade $Q_{{\bold 
F},n}(t).$
Then $Q_n(t) := Q_{\beta,n}(t)Q_{{\bold F},n}(t)$ is a 
weighted cascade,
with weights ${\bold F}_\beta$ defined by $(F_\beta)_\gamma
:= B_\gamma F_\gamma,$  such that for $\sigma\in M^+(T)$
$$Q_\infty\sigma  = Q_{\beta,\infty}( Q_{{\bold 
F},\infty}\sigma)\quad
  a.s.$$
\endproclaim

One way in which the significance of this
approach is illustrated is by  the existence of the 
following
 natural weight systems.

\proclaim{Theorem 2.5 \rm(Kahane decomposition)} The 
Kahane decomposition
\cite{K2, K4}, \break
$Q_n(t) = Q_n^\prime(t) + Q_n^{\prime\prime}(t),$ is 
equivalent to
a weight decomposition with
respect to the weights defined by $F_{t|n} = 
{Q_n^\prime(t)\over Q_n(t)}, t\in
T,$  and $1 - F_{t|n}$ for $Q_n^\prime(t), 
Q_n^{\prime\prime}(t)$
respectively.
 \endproclaim

\proclaim {Theorem 2.6 \rm(Dimension decomposition)} The
 decomposition \cite{KK, C} $\lambda_\infty = \int
(\lambda_\infty)_{\alpha}\nu\,(d\alpha)$ is equivalent to 
a weight decomposition
with respect to the weight systems ${\bold F}_s := 
\{F_{s,\gamma}:\gamma\in
T^*\}, s\in [0,1],$ defined by 
$$F_{s,\gamma} = {E[\nu_\gamma([0,s]) | {\Cal F}_n]\over
\lambda_n(\Delta_\gamma)},$$ where
$$\nu_\gamma([0,s]) :=
 \int_{[0,s]}(\lambda_\infty)_\alpha(\Delta_\gamma)\nu\,(d%
\alpha),
\quad \quad n = |\gamma|.$$
\endproclaim

\rem{Remark \rm 1} Two interesting classes of dependent 
cascades 
are: (i) {\it Markov cascades and,} (ii) {\it exchangeable 
cascades.} 
 One may compute the weighting decompositions 
for  Markov cascades, for example, in terms
of {\it harmonic measures} corresponding to hitting 
probabilities of 
the {\it survival classes} introduced in \cite{WW1}.
 Also for  exchangeable cascades  one may compute these 
weighting
 decompositions in terms of  conditioned de Finetti 
measures applied to
 a partition of $M^+(T).$\endrem

\rem{Remark \rm 2}  One may show that the maps
 $ s\to F_{s,\gamma}$
are nondecreasing in $s$ and the Lebesgue-Stieltjes 
measures  associated
 with the maps
$s\to \int F_{s,t|n}\lambda_n\,(dt)$ converge vaguely to 
the spectral measure
$\nu. $ 
In fact, if  the map
$ s\to F_{s,\gamma}$ is absolutely continuous with respect 
to some measure
$\sigma,$ then $t\to \int{dF_{t,\gamma}\over 
d\sigma}\lambda_n\,(d\gamma )$ is
itself a positive $T=[0,\infty)$-martingale. In this case
the percolation method of Theorem 2.4 may be applied to 
compute the
dimension spectrum. This, in fact, illustrates a more 
general notion
of {\it differentiable weights} which are introduced  to 
compute more 
detailed 
dimension spectra in \cite{WW3}.

Finally, the scope of the theory is well illustrated by 
the following
approach to the Kahane-Peyri\`ere theorem.  In particular, 
the 
suitability of this method for more general cascades is 
made transparent by
the  suppressed role  of independence beyond the 
martingale structure and/or
the general ergodicity under the size-biased distribution.  

Consider the independent cascade. Let $W$ denote a generic 
generator
with distribution $q(dx).$
We write ${\bold F} = {\bold 1}$ to denote the case of 
unit weights.\endrem

\demo \nofrills {\bf A sufficient condition for 
nondegeneracy.}\quad
Let us first show how it follows that
  $E_P W\log_b W < 1$ is a sufficient condition for 
nondegeneracy. 
 Choose 
$0 < \log_b c < 1 -  E_P W\log_b W.$  Fix $t\in T.$ 
Under $P_{{\bold 1},t},$
for $\gamma$\<'s  along the $t$-path,
i.e. $\gamma = t|j,$ some j,
 the  $ W_\gamma $\<'s  
are i.i.d.\ with distribution $xq(dx);$ thus the name 
\lq\lq sized-biased\rq\rq.
 For $\gamma$\<'s off
of\ the $t$-path, i.e. 
 $\gamma \neq t|j,$ the
$W_\gamma$\<'s are i.i.d.\ with distribution $q(dx).$  It 
follows from
the SLLN (i.e. ergodicity) that $P_{{\bold 1},t}$-a.s.,  
$\root j \of {{c^j\over b^j}\prod_{i=1}^jW_{t|i}}
\to {c\over b}e^{E_PW\log W}.$  Now use Theorem 2.3,
with $c_j = c^{-j},$  to conclude
that $\lambda_\infty(T) < \infty, P_{{\bold 1},t}$-a.s., and
therefore, by Theorem 2.2, after integrating out t, 
$\lambda_\infty(T) < \infty, {\Cal Q}_{{\bold 1}}$-a.s. 
Now apply Theorem 2.1 to  see  that $E\lambda_\infty(T) = 
1.$
\enddemo
\demo \nofrills {\bf A necessary  condition for 
nondegeneracy.}\quad
To  see why  the entropy condition $E_P W\log_b W < 1$ is 
also necessary
for nondegeneracy,
suppose that $E_P W\log_b W  \ge  1, P_t(W = b) < 1.$
Fix $t\in T.$ If $E_P W\log_b W >  1$, then one obtains
$\lambda_n(T) \to \infty$  $P_{{\bold 1},t}$-a.s.  
directly from Theorem 2.3.
If $ E_P W\log_b W = 1$ but $ P_t(W = b) < 1$,  then, 
since under
the size-biasing $W > 0$ a.s., one has that
$\limsup\log\lambda_n(T) = \infty,  P_{{\bold 1},t}$-a.s.\
 by the Chung-Fuchs null-recurrence
criterion for the random walk along the $t$-path. Integrate
out $t$ and use Theorems 2.2 and 2.1 (in that order) to 
conclude
that $\lambda_\infty(T) = 0, P$-a.s.  The special case of 
the degeneracy
 $P_t(W=b) = 1,$ i.e.\ $ P(W=b)={1\over b}=1-P(W=0),$  may
be handled similarly or directly from the theory of 
branching
processes. 
\enddemo
In a similar manner one may obtain the Kahane-Peyri\`ere 
divergence
of moments criterion from Theorems 2.3, 2.2, and 2.1 (in 
that order).
The condition for nondegeneracy and the percolation method 
of  Theorem
2.4 may then be combined to obtain the Hausdorff dimension 
of 
the support of $\lambda_\infty$  in the manner  first 
noted in [K1].

The detailed proofs for the underlying theory and 
applications to
general classes of dependent cascades  appear in  the  
companion
paper  \cite{WW3}.

\enddocument